\documentclass{amsart}

\newcommand{\ignore}[1]{}

\usepackage{graphicx}

\newcommand{\enum}[1]{\begin{enumerate} {#1} \end{enumerate}}

\newcommand{\dfn}[1]{\textbf{#1}}
\newcommand{\xquote}[1]{\begin{quote} #1 \end{quote}} %x b/c \quote already defined%

\newcommand{\figref}{Figure~\ref}

\newcommand{\st}{:}
\newcommand{\x}{\times}
\newcommand{\nsub}{immediate subset}
\newcommand{\nsup}{immediate superset}
\newcommand{\isomorphic}{\simeq}
\newcommand{\setsubtr}{\backslash}  % set subtraction

\newtheorem{theorem}{Theorem}
\newcommand{\thm}[1]{ \begin{theorem} {#1} \end{theorem} }
\newtheorem{corollary}[theorem]{Corollary}

\newtheorem{lemma}[theorem]{Lemma}
\newcommand{\lem}[1]{ \begin{lemma} {#1} \end{lemma} }
\newtheorem{proposition}[theorem]{Proposition}
\newcommand{\prop}[1]{ \begin{proposition} {#1} \end{proposition} }

\theoremstyle{definition} %%{
\newtheorem{definition}{Definition}

\newtheorem{example}{Example}

\theoremstyle{plain} %%}

\newcommand{\pf}[1]{\begin{proof}{#1}\end{proof}}

\begin{document}{

\title{Induced Subgraphs of Johnson Graphs}
\author{Ramin Naimi}
\address{
Department of Mathematics,
Occidental College,
Los Angeles, CA 90041, USA.
}
\author {Jeffrey Shaw}
\address{
Department of Mathematics,
Occidental College,
Los Angeles, CA 90041, USA.
}
\date{\today}

\begin{abstract}
The Johnson graph $J(n,N)$ is defined as
the graph whose vertices are
the $n$-subsets of the set $\{1, 2, \cdots, N\}$,
where two vertices are adjacent
if they share exactly $n-1$ elements.
Unlike Johnson graphs,
induced subgraphs of Johnson graphs (JIS for short)
do not seem to have been studied before.
We give some necessary conditions and some sufficient conditions
for a graph to be JIS, including:
in a JIS graph, any two maximal cliques share at most two vertices;
all trees, cycles, and complete graphs are JIS;
disjoint unions and Cartesian products of JIS graphs are JIS;
every JIS graph of order $n$
is an induced subgraph of $J(m,2n)$
for some $m \le n$.
This last result gives an algorithm for deciding if a graph is JIS.
We also show that all JIS graphs are
edge move distance graphs, but not vice versa.
\end{abstract}

\maketitle

%%%%%  Statements of theorems that are listed more than once (in Intro and other sections)

\newcommand{\CliquesInJISTheoremStatement}
{Suppose $G$ is JIS and
$L$ and $L'$ are distinct maxcliques in $G$.
Then:

\enum{

\item
$L$ and $L'$
share at most two vertices.

\item
If $L$ and $L'$ share exactly two vertices, then no vertex in $V(L)
\setsubtr V(L')$ is adjacent to a vertex in $V(L') \setsubtr V(L)$.

\item
If $L$ and $L'$ share exactly one vertex,
then each vertex in either of the two sets
$V(L) \setsubtr V(L')$ and $V(L') \setsubtr V(L)$
is adjacent to at most one vertex in the other set.

}
}

\newcommand{\CliqueChainTheoremStatement}
{Suppose  $L_1,  \cdots, L_k$, where $k$ is odd and at least 3,
are distinct maxcliques
in a graph $G$
such that
$L_i$ shares exactly two vertices with $L_{i+1}$ for $1 \le i \le k-1$,
and $L_k$ shares exactly two vertices with $L_1$;
then $G$ is not JIS.
}

\newcommand{\CompletesAreJISTheoremStatement}
{All complete graphs are JIS.}

\newcommand{\CompletesAndCyclesAreJISTheoremStatement}
{All complete graphs and all cycles are JIS.}

\newcommand{\CyclesAreJISTheoremStatement}
{All cycles are JIS.}

\newcommand{\IffCoreIsJISTheoremStatement}
{A graph is JIS iff its 2-core is JIS.}

\newcommand{\ConnectedComponentsTheoremStatement}
{A graph is JIS iff all its connected components are JIS.}

\newcommand{\CartesianProductsTheoremStatement}
{The Cartesian product of two JIS graphs is JIS.}

\newcommand{\UpperboundTheoremStatement}
{Every JIS graph of order $n$
is isomorphic,
for some $m \le n$,
to an induced subgraph of the Johnson graph $J(m,2n)$.
}

\newcommand{\JISimpliesEMDGTheoremStatement}
{Every JIS graph is an edge move distance graph.}

\newcommand{\KnMinusOneEdgeIsNotJISTheoremStatement}
{The graph obtained by
removing one edge from the complete graph $K_n$,
where $n \ge 5$,
is an edge move distance graph but is not JIS.}

\section{Introduction}

We work with finite, simple graphs.
Let $F = \{S_1, \cdots, S_m\}$
be a family of finite sets.
The \dfn{intersection graph} of $F$,
denoted $\Omega(F)$,
is the graph whose vertices are the elements of $F$,
where two vertices  $S_i$ and $S_j$, $i \ne j$,
are adjacent if they share at least one element.
More generally,
for a fixed positive integer $p$,
the \dfn{$p$-intersection graph} of $F$,
denoted $\Omega_p(F)$,
is the graph whose vertices are the elements of $F$,
where two vertices
are adjacent if they share at least $p$ elements.
(Thus $\Omega_p(F)$ is a subgraph of $\Omega_1(F) = \Omega(F)$.)
McKee and McMorris \cite{MKMM}
give an extensive and excellent survery of intersection graphs,
which also includes a section on $p$-intersection graphs.
Here we narrow attention to
$p$-intersection graphs of families of $(p+1)$-sets,
so that
two vertices $S_i$ and $S_j$ are adjacent
if $|S_i \cap S_j| = |S_i|-1 = |S_j|-1$,
i.e., $S_i$ and $S_j$
differ by exactly one element.

Another way to view these graphs is as induced subgraphs of Johnson graphs.
Given positive natural numbers $n \le N$,
the \dfn{Johnson graph $J(n,N)$} is defined as
the graph whose vertices are
the $n$-subsets of the set $\{1, 2, \cdots, N\}$,
where two vertices are adjacent
if they share exactly $n-1$ elements.
Hence a graph G is isomorphic to
an induced subgraph of a Johnson graph
iff it is possible to assign,
for some fixed $n$,
an $n$-set $S_v$
to each vertex $v$ of $G$ such that
distinct vertices have distinct corresponding sets,
and vertices $v$ and $w$ are adjacent iff
$S_v$ and $S_w$ share exactly $n-1$ elements.
When this happens,
we say the family of $n$-sets $F=\{S_v \st v \in V(G) \}$
\dfn{realizes} $G$ as
an induced subgraph of a Johnson graph,
which we abbreviate by saying $G$ is \dfn{JIS}.
Thus, $F$ realizes $G$ as a JIS graph iff
$G$ is isomorphic to $\Omega_{n-1}(F)$,
which in turn is isomorphic to an induced subgraph of $J(n,N)$,
where $N= | \bigcup_{S \in F} S |$.

Although there is a considerable amount of literature
written on Johnson graphs,
we have not been able to find any
on their induced subgraphs.
It would be desirable to obtain
``nice'' necessary and sufficient conditions
for when a graph is JIS.
In this paper, we only give
some necessary conditions and some sufficient conditions.

A \dfn{clique} in a graph $G$ is a complete subgraph of $G$.
A clique $L$ in $G$ is called a
\dfn{maximal clique}, or a \dfn{maxclique} for short,
if there is no larger clique $L' \subseteq G$
that contains $L$.
In Section~\ref {SecMaxcliques}
we describe how the maxcliques of a graph play a role in
whether or not it is JIS.
In particular, part~(1) of Proposition~\ref{CliquesInJISTheorem}
states that any two distinct maxcliques in a JIS graph can share at most two vertices.
It follows, for example, that
the graph ``$K_5$ minus one edge''
is not JIS,
since it contains two maximal 4-cliques that share three vertices.

The conditions given in Section~\ref{SecMaxcliques}
are necessary, but not sufficient,
for a graph to be JIS.
In Section~\ref{SecMiscJISgraphs}
we show that
the complete bipartite graph $K_{2,3}$,
as well as a few other graphs,
satisfy all these necessary conditions
but are not JIS.
In Section~\ref{SecMiscJISgraphs}
we also give some sufficient conditions for a graph to be JIS,
including the following:
\CompletesAndCyclesAreJISTheoremStatement \
\ConnectedComponentsTheoremStatement \
\CartesianProductsTheoremStatement \

Despite not having a ``nice'' characterization of JIS graphs,
for any graph $G$ the question ``Is $G$ JIS?'' is decidable;
this follows from Theorem~\ref{UpperboundTheorem},
which says:
\UpperboundTheoremStatement \
In other words,
every JIS graph of order $n$
can,
for some $m \leq n$,
be realized by $m$-subsets of
$\{1, 2, \cdots, 2n\}$.
This gives us a simple (albeit slow) algorithm
for determining if a graph $G$ is JIS:
Do an exhaustive search among all $n$-families of
$m$-subsets of $\{1, \cdots, 2n\}$,
where $n$ is the order of $G$ and $m \le n$,
to see if any of them realizes $G$ as a JIS graph.

The \emph{$p$-intersection number} of a graph $G$ is defined as
the smallest $k$ such that
$G$ is isomorphic to the $p$-intersection graph of
a family of subsets of $\{1, \cdots, k\}$ (\cite{MKMM}, p.~91).
Thus, an immediate corollary of Theorem~\ref{UpperboundTheorem} is that
every JIS graph of order $n$ has,
for some $m \le n$,
$(m-1)$-intersection number at most $2n$.

In the final section of this paper
we discuss edge move distance graphs
and their relationship to JIS graphs.

\section{Maxcliques in JIS Graphs}
\label{SecMaxcliques}

Given $n$-sets $S_1, \cdots, S_k$
with $n \ge 1$ and $k \ge 2$,
we say they share an
\dfn{\nsub}
if  $|\bigcap_{i=1}^k S_i | = n-1$.
Similarly,
$S_1, \cdots, S_k$ share an
\dfn{\nsup}
if  $|\bigcup_{i=1}^k S_i | = n+1$.
Observe that for $k=2$,
$S_1$ and $S_2$ share an \nsub\
iff they share an \nsup:
$|S_1 \cup S_2| = |S_1| + |S_2| - |S_1 \cap S_2| = 2n - |S_1 \cap S_2|$;
hence
$|S_1 \cup S_2| = n+1$ iff
$|S_1 \cap S_2| = n-1$.
We begin with the following elementary result
on realizations of complete graphs as JIS graphs.

\lem{ \label{lem1}
Let $S_1, \cdots, S_k$ be $n$-sets that
pairwise share an \nsub,
where $n \ge 1$ and $k \ge 3$.
Then $S_1, \cdots,  S_k$ share
an \nsub\
or an \nsup,
but not both.
}

\pf{

We first show that for $k \ge 3$,
if $S_1, \cdots, S_k$
share an \nsub,
then they do not share an \nsup.
Suppose $T = S_1 \cap \cdots \cap S_k$
has $n-1$ elements.
Then, for each $i$,
$S_i \setsubtr T$ has exactly one element, $x_i$.
For all $j \ne i$,
$x_i \not\in S_j$ since $S_i \ne S_j$.
It follows that
$S_1 \cup \cdots \cup S_k$ has at least
$n-1+k \ge n+2$ elements,
since $k \ge 3$.
Thus
$S_1, \cdots, S_{k}$ do not share
an \nsup.

Now suppose $S_1, \cdots, S_{k}$
pairwise share an \nsub.
We use induction on $k$ to prove
that they share an \nsub\ or an \nsup.

Assume $k=3$.
Let $T = S_1 \cap S_2$.
If $T \subset S_3$,
then $|S_1 \cap S_2 \cap S_3| = |T| = n-1$,
and we're done.
So assume $T \not \subset S_3$.
Note that
$|S_1 \setsubtr T| = |S_2 \setsubtr T| = 1$.
Hence, for $S_3$ to share $n-1$ elements
with each of $S_1$ and $S_2$,
it must contain an $(n-2)$-subset of $T$,
as well as $S_1 \setsubtr T$ and $S_2 \setsubtr T$,
and no other elements.
It follows that
$|S_1 \cup S_2 \cup S_3| = n+1$,
as desired.

Now assume $k\ge4$.
Then, by our induction hypothesis,
$S_1, \cdots, S_{k-1}$ share
an \nsub\ or an \nsup;
and similarly for $S_2, \cdots, S_{k}$.
We have four cases:

Case 1:
$S_1, \cdots, S_{k-1}$ share an \nsub\
and
$S_2, \cdots, S_{k}$ share an \nsub.
Then
$S_1, \cdots, S_{k}$ share $S_2 \cap S_3$
as an \nsub.

Case 2:
$S_1, \cdots, S_{k-1}$ share an \nsup\
and
$S_2, \cdots, S_{k}$ share an \nsup.
Then
$S_1, \cdots, S_{k}$ share $S_2 \cup S_3$
as an \nsup.

Case 3:
$S_1, \cdots, S_{k-1}$ share an \nsub\
and
$S_2, \cdots, S_{k}$ share an \nsup.
Let $T = S_1 \cap \cdots \cap S_{k-1}$.
Then,
for $1 \le i \le k-1$,
$S_i \setsubtr T$ has exactly one element, $x_i$;
and, for  $1 \le j \le k-1$ with $j \ne i$,
$x_i \not\in S_j$ since $S_i \ne S_j$.
Since $|S_2 \cup \cdots \cup S_k| = n+1 = |S_2 \cup S_3|$,
$S_k$ is a proper subset of
$S_2 \cup S_3 = T \cup \{x_2, x_3\}$.
And since $S_2$, $S_3$, $S_k$
share an \nsup,
they do not share an \nsub;
hence $T \not \subset S_k$.
This implies that $x_2, x_3 \in S_k$
since $S_k$ has $n$ elements
and $T \cup \{x_2, x_3\}$ has $n+1$ elements.
But $x_2, x_3 \not \in S_1$,
so $|S_1 \cap S_k| < n-1$,
which contradicts the hypothesis of the lemma.

Case 4:
$S_1, \cdots, S_{k-1}$ share an \nsup\
and
$S_2, \cdots, S_{k}$ share an \nsub.
This case is similar to Case~3.

}

We now use Lemma~\ref{lem1}
to establish restrictions on
how maxcliques in a JIS graph
can intersect or connect to each other by edges.

\prop{ \label{CliquesInJISTheorem}
\CliquesInJISTheoremStatement
}

\pf{
Let $\{S_v \st v \in V(G)\}$
be a family of $n$-sets that
realizes $G$ as a JIS graph.

\medskip
\emph{Proof of Part 1.}
Suppose towards contradiction that
$L$ and $L'$ are distinct maxcliques
that share three (or more) vertices, $u$, $v$, and $w$.
Let $x$ be a vertex of $L$ not in $L'$,
and $x'$ a vertex of $L'$ not in $L$;
$x$ and $x'$ exist since $L$ and $L'$
are distinct and maximal.
Then, by Lemma~\ref{lem1},
the sets $S_x$, $S_u$, $S_v$, and $S_w$
share an \nsub\ or an \nsup.
Similarly for  $S_{x'}$, $S_u$, $S_v$, and $S_w$.
But $S_u$, $S_v$, and $S_w$ cannot share
both an \nsub\ and an \nsup.
It follows that
$S_x$ and $S_{x'}$
share an \nsub\ or an \nsup,
which implies that $x$ and $x'$ are adjacent.
Hence every vertex of $L$
is adjacent to every vertex of $L'$,
but this contradicts the assumption that $L$
is a maxclique in $G$.

\medskip
\emph{Proof of Part 2.}
Let $L$ and $L'$ be distinct maxcliques
that share exactly two vertices, $v$ and $w$.
Suppose towards contradiction that
there exist adjacent vertices
$x \in V(L) \setsubtr V(L')$
and
$x' \in V(L') \setsubtr V(L)$.
Then the induced subgraph of $G$ containing
$\{x, x', v, w\}$ is a 4-clique.
Let $L''$ be the maxclique that contains
this 4-clique.
Then $L''$ is distinct from $L$
and shares at least three vertices with it.
This contradicts Part~1.

\medskip
\emph{Proof of Part 3.}
The proof is similar to the proof of Part~2.
Let $L$ and $L'$ be distinct maxcliques
that share exactly one vertex, $v$.
Suppose towards contradiction that
there exist  vertices
$x \in V(L) \setsubtr V(L')$
and
$x', y'  \in V(L') \setsubtr V(L)$
with $x$ adjacent to $x'$ and $y'$.
Then the induced subgraph of $G$ containing
$\{x, x', y', v\}$ is a 4-clique,
and the maxclique that contains
this 4-clique
is distinct from $L'$
and shares at least three vertices with it.
This contradicts Part~1.

}

\prop{ \label{CliqueChainTheorem}
\CliqueChainTheoremStatement
 }

\pf{
In the following,
$L_{i+1}$ refers to $L_1$ whenever $i=k$.
Suppose towards contradiction that
$G$ is realized as a JIS graph
by a family of $n$-sets.
Note that each $L_i$ has at least three vertices,
since otherwise it would not be distinct from $L_{i+1}$.
Hence, by Lemma~\ref{lem1},
we can label each $L_i$ as
either ``sub'' or ``super''
according to whether
the $n$-sets assigned to its vertices
share an \nsub\ or an \nsup.
Then, since $k$ is odd,
there exists a $j$ such that
$L_j$ and $L_{j+1}$ have the same label.
Now, $L_j$ and $L_{j+1}$ share two vertices;
therefore the $n$-sets assigned to their vertices
must all share the same \nsub\ or \nsup,
which makes all vertices in $L_j$ adjacent
to those in $L_{j+1}$,
giving a contradiction.

}

An equivalent way of stating the above result is:
One can label every maxclique in a JIS graph
with a ``$+$'' or ``$-$'' (or any two symbols)
in such a way that
any two maxcliques that share two vertices
have distinct labels.

\section{Miscellaneous JIS and non-JIS graphs}
\label{SecMiscJISgraphs}

In this section
we give some sufficient conditions
for when a graph is JIS.
We also describe some graphs
that satisfy all the conditions
listed in the results of the previous section
as necessary for a graph to be JIS,
but are not JIS.

\prop{ \label{CompletesAreJISTheorem}
\CompletesAndCyclesAreJISTheoremStatement
}

\pf{
For each $n$,
$K_n$ is realized as a JIS graph by the 1-sets
$\{1\}$, $\{2\}$,  $\cdots$, $\{n\}$.
For each $n \ge 3$,
the $n$-cycle is realized as a JIS graph by the 2-sets
$\{1,2\}$, $\{2,3\}$, $\cdots$, $\{n-1, n\}$,  $\{n, 1\}$.
}

We define the \dfn{$n$-core} of a graph $G$
as the graph obtained by recursively
removing all vertices of degree less than $n$
until there are none left.

\prop{ \label{IffCoreIsJISTheorem}
\IffCoreIsJISTheoremStatement
}

\pf{
Suppose $G$ is obtained from a graph $G'$ by removing
exactly one vertex, $w$,
which has degree~0 or 1.
By induction, it is enough to show that
$G$ is JIS iff $G'$ is JIS.
Clearly, if $G'$ is JIS,
then so is $G$,
since any induced subgraph of a JIS graph is JIS.
To prove the converse,
suppose $G$ is JIS.
Let $\{ S_x \st x \in V(G) \}$
be $n$-sets that realize $G$ as a JIS graph.
Pick distinct $a$ and $b$
that are not in any of the sets $S_x$.
For each $x \in V(G)$,
let $S'_{x} = S_x \cup \{a\}$.
Let $S'_w = S_v \cup \{b\}$,
where $v \in V(G')$ is arbitrary
if $w$ has degree 0,
and $v$ is adjacent to $w$
if $w$ has degree 1.
Then $\{ S'_x \st x \in V(G') \}$
are $(n+1)$-sets that realize $G'$
as a JIS graph,
as desired.

}

It follows as a trivial corollary that all trees are JIS.

\prop{ \label{ConnectedComponentsTheorem}
\ConnectedComponentsTheoremStatement
}

\pf{
One direction is trivial:
every induced subgraph of a JIS graph,
and in particular every connected component of it,
is JIS.
We prove the converse by induction on the number of components of $G$.

Base step:
Suppose that $G$ has two components,
$G_i$, $i=1, 2$,
each realized as a JIS graph
by a family of sets $F_i$.
We can assume without loss of generality that
each set in $F_1$ is disjoint
from each set in $F_2$.

We would like each set in $F_1$
to have the same size as
each set in $F_2$,
in order to obtain $F_1 \cup F_2$
as a family that realizes $G$ as a JIS graph.
If this is not already so, we proceed as follows.
Let $m_i$ denote the number of elements in each set in $F_i$.
We can assume $n_1 > n_2$.
Now add the first $n_1 - n_2$ elements of
the first set in $F_1$
to every set in $F_2$.

Once the sets in the two families all have the same size,
we must make sure that
sets corresponding to vertices in different components of $G$
do not share immediate subsets.
This will automatically be true
for sets that had two or more elements
before any extra elements were added to them
(since we started with the sets in $F_1$ disjoint from those in $F_2$),
but not for singletons.
We remedy this
by adding, for each $i$,
an element $e_i$
to every set in $F_i$,
where $e_1$ and $e_2$
are distinct elements not already in any set in any $F_i$.
It is now easy to verify that
$F_1 \cup F_2$ realizes $G$ as a JIS graph.

The inductive step follows trivially from the base step.
}

\prop{ \label{CartesianProductsTheorem}
\CartesianProductsTheoremStatement
}

\pf{
Let $G$ and $G'$ be JIS graphs
that are realized, respectively,
by sets $\{ S_x \st x \in V(G) \}$ and
$\{ S'_{x'} \st x' \in V(G') \}$.
We can assume without loss of generality
that every $S_x$ is disjoint from every $S'_{x'}$.

For each vertex $v=(x, x') \in V(G \x G')$,
let $T_v = S_x \cup S'_{x'}$.
By definition, two vertices
$v=(x, x')$  and $w=(y, y')$ of $G \x G'$
are adjacent iff
$x=x'$ and $y$ is adjacent to $y'$
or
$y=y'$ and $x$ is adjacent to $x'$.
Thus, $T_v$ and $T_w$ share an \nsub\ iff
$v$ and $w$ are adjacent.
Hence the sets $\{T_v \st v \in G \x G' \}$
realize $G \x G'$ as a JIS graph.
}

\prop{ \label{K23Theorem}
The complete bipartite graph
$K_{2,3}$ is not JIS.
}

\begin{figure}[ht]

 \centering
 \includegraphics[width=50mm]{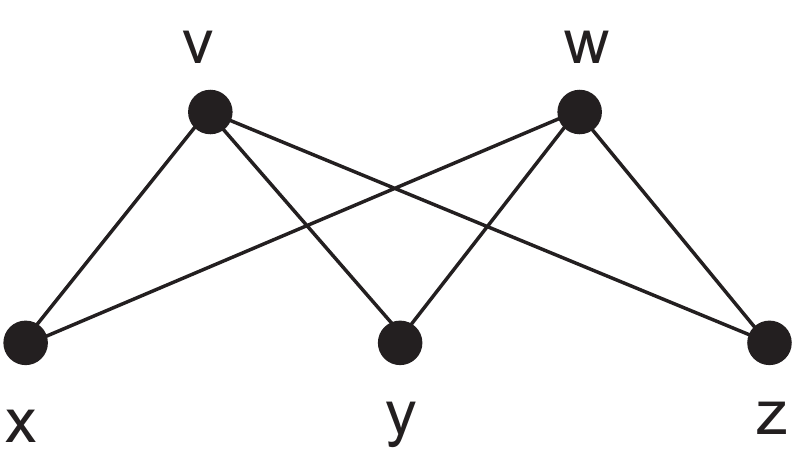}
 \caption{$K_{2,3}$ with labeled vertices
 \label{KTwoThreeLabeledFigure}
 }
\end{figure}

\pf{
Label the two degree-3 vertices of $K_{2,3}$ as $v$ and $w$,
and the three degree-2 vertices as $x$, $y$, and $z$,
as in \figref{KTwoThreeLabeledFigure}.
Suppose towards contradiction that
there exists a family of $n$-sets
$\{S_u \st u \in V(K_{2,3})\}$
that realizes $K_{2,3}$ as a JIS graph.
Since $v$ and $w$ have distance two
(where \emph{distance} is the number of edges
in the shortest path joining the two vertices),
$S_v$ and $S_w$ must share exactly $n-2$ elements
(this does not work for distance $\ge 3$;
it works only for distance $\le 2$).
Let $T = S_v \cap S_w$.
Then, since each of $x$, $y$, and $z$
is adjacent to both $v$ and $w$,
$S_x$, $S_y$, and $S_z$ must each contain $T$ as a subset.
Therefore, by subtracting $T$ from every $S_u, u \in V(K_{2,3})$,
we get a family of 2-sets that realizes $K_{2,3}$.
Hence we will assume that every $S_u$ has exactly two elements.
It follows that $S_v$ and $S_w$ are disjoint;
and $S_x$, $S_y$, and $S_z$
are pairwise disjoint
and each shares
exactly one element with each of $S_v$ and $S_w$.

So, without loss of generality,
$S_v = \{1,2\}$, and
$S_w = \{3,4\}$.
Therefore, again without loss of generality,
$S_x = \{1,3\}$, and
$S_y = \{2,4\}$.
And there is nothing left for $S_z$.
}

The graph $K_{2,3}$ can be thought of as
two 4-cycles that share three vertices.
So one may wonder whether
the graph $\theta_n$ consisting of
two $n$-cycles that share $n-1$ vertices
is also not JIS.
It turns out that $\theta_n$ is not JIS
only for $n=4$ and $n=5$.
The proof that $\theta_5$ is not JIS
is very similar to the proof that
$K_{2,3}$ is not JIS,
and we therefore omit it.
The proof that $\theta_n$ is JIS for  $n \ge 6$
is a straightforward construction, which we also omit.

One may also wonder whether $K_{2,3}$ becomes JIS
if an edge is added to it.
There are, up to isomorphism,
two ways to add an edge to $K_{2,3}$:
add an edge that connects the two degree-3 vertices; or
add an edge that connects two of the three degree-2 vertices.
It turns out that neither of these two graphs is JIS.
The proof that the former graph is not JIS
follows immediately from
Proposition~\ref{CliqueChainTheorem}.
The proof that the latter graph
(which we call $\Delta_2$)
is not JIS is given below
in Proposition~\ref{Delta2and4NotJISTheorem}.

The graphs $\Delta_i$ depicted in \figref{Delta234Figure}
have the following pattern
(ignore the vertex labels and the $+$ and $-$ signs for now;
they are used later):
$\Delta_i$ consists of a chain of $i$ ``consecutively adjacent'' triangles,
plus one vertex which is connected to
the two vertices of degree~2 in the triangle chain.
It turns out that,
like $K_{2,3}$,
$\Delta_2$, $\Delta_4$, and $\Delta_6$
satisfy the necessary conditions
in the results of the previous sections
for being JIS,
but are not JIS;
$\Delta_3$ and $\Delta_5$, however,
are JIS.
We prove these claims below,
except for $\Delta_6$:
its proof is similar to that of $\Delta_2$ and $\Delta_4$,
but is more tedious,
and in our opinion not worth being included here.
We did not check which $\Delta_i$ are JIS for $i \ge 7$,
but, from the pattern for $i \le 6$, it seems that:

\xquote{
\emph{Conjecture:}
$\Delta_i$ is JIS iff $i$ is odd.
}

\begin{figure}[ht]

 \centering
 \includegraphics[width=100mm]{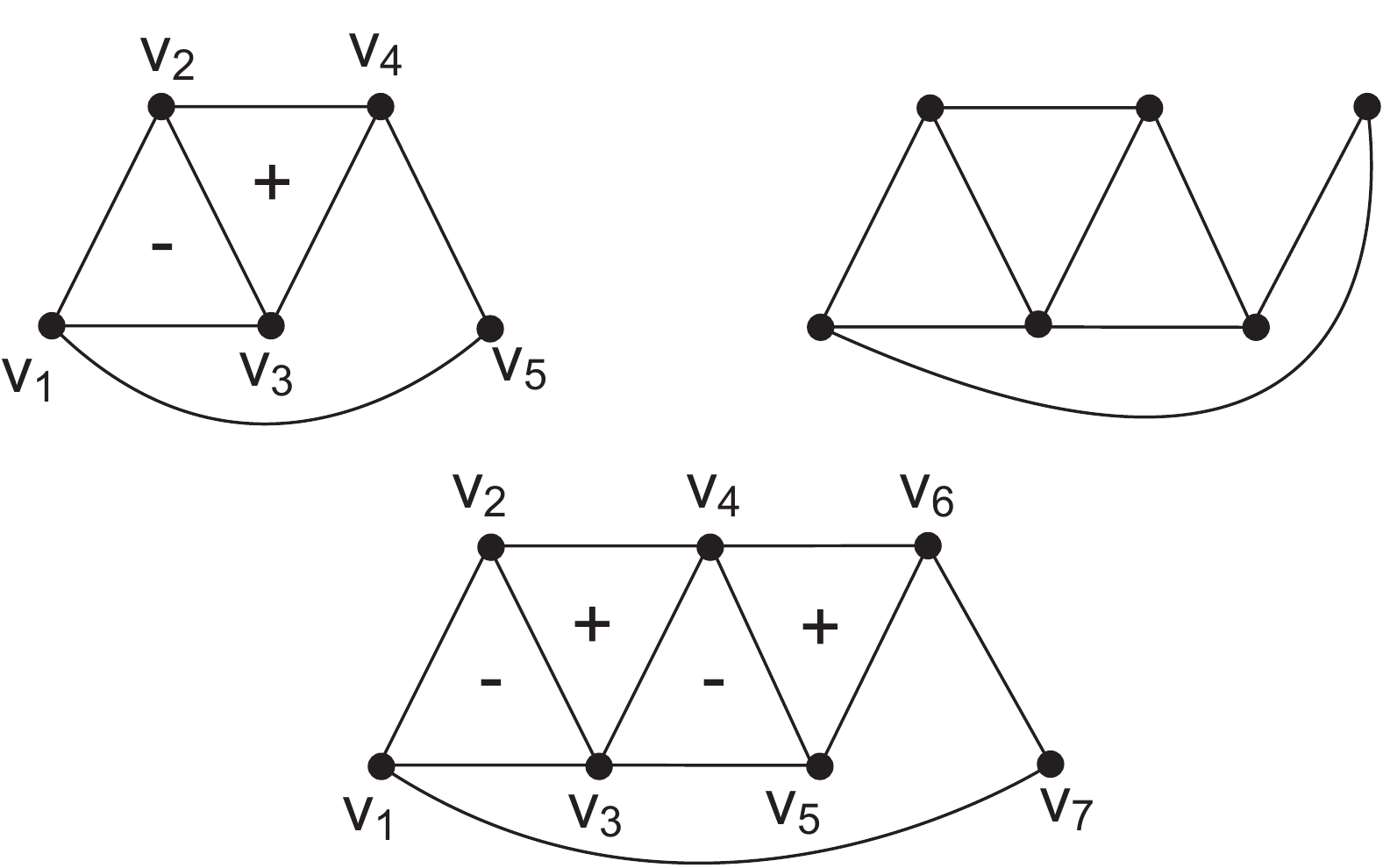}
 \caption{$\Delta_2$, $\Delta_3$, and $\Delta_4$,
 with vertices labeled in $\Delta_2$ and $\Delta_4$.
 \label{Delta234Figure}
 }
\end{figure}

\prop{
\label{Delta2and4NotJISTheorem}
(i) The graphs $\Delta_2$ and $\Delta_4$ are not JIS.
(ii) The graphs $\Delta_3$ and $\Delta_5$ are JIS.
}

Remark:
As mentioned above,
$\Delta_2$ is isomorphic to $K_{2,3}$ plus an edge
that connects two of its three degree-2 vertices.
Because of this,
the proof that $K_{2,3}$ is not JIS
can be easily modified to prove that
$\Delta_2$ is not JIS.
However, we give a different proof below,
one that can be naturally extended
to also prove that
$\Delta_4$ (and $\Delta_6$) is not JIS.

\pf{
Label the vertices of $\Delta_2$
as $v_1, \cdots, v_5$, as in \figref{Delta234Figure}.
The $+$ and $-$ signs will be explained shortly.
Suppose, towards contradiction,
that $\Delta_2$ can be realized as a JIS graph
by sets $S_1, \cdots, S_5$
(for simplicity, we write $S_i$ instead of $S_{v_i}$).
Each of the two triangles in $\Delta_2$
is a maxclique.
Thus, by Lemma~\ref{lem1},
$S_1$, $S_2$, and $S_3$ must share an \nsub\ or an \nsup;
similarly for $S_2$, $S_3$, and $S_4$.
Furthermore,
$S_1$, $S_2$, and $S_3$ share an \nsub\ iff
$S_2$, $S_3$, and $S_4$ share an \nsup,
because:
if
$S_1$, $S_2$, and $S_3$ share an \nsub\ and
$S_2$, $S_3$, and $S_4$ also share an \nsub,
then
$S_1$ and $S_4$ must share
$S_2 \cap S_3$ as an \nsub,
but this contradicts the fact that
$v_1$ and $v_4$ are not adjacent;
and if
$S_1$, $S_2$, and $S_3$ share an \nsup\ and
$S_2$, $S_3$, and $S_4$ also share an \nsup,
then
$S_1$ and $S_4$ must share
$S_2 \cup S_3$ as an \nsup,
which implies that they also share an \nsub,
again contradicting the fact that
$v_1$ and $v_4$ are not adjacent.

Thus, without loss of generality,
we will assume that
$S_1$, $S_2$, and $S_3$ share an \nsub.
This is indicated in \figref{Delta234Figure} by the $-$ sign;
the $+$ signs indicate \nsup s.
So we will assume that
$S_1 = \{1,2,3,4\}$,
$S_2 = \{1,2,3,5\}$, and
$S_3 = \{1,2,3,6\}$;
we explain in the next paragraph why
there is no loss of generality in assuming that
$S_i$ are 4-sets (as opposed to larger sets).
To make the notation more compact,
we will drop the commas and the braces from each set;
e.g., $S_1 = 1234$.
Then $S_4$ must be a 4-subset of
$S_2 \cup S_3 = 12356$.
Since $S_1$ and $S_4$ have no \nsub,
we can without loss of generality assume that
$S_4 = 2356$.
Now, $S_5$ must differ by exactly one element
from each of $S_1$ and $S_4$.
The only possibilities are
1235, 1236, 2345, and 2346.
But the first two are equal to $S_2$ and $S_3$ respectively;
and the last two differ from $S_2$ and $S_3$ respectively
by exactly one element,
which is not allowed
since $v_5$ is adjacent to neither $v_2$ nor $v_3$.
Thus we have a contradiction, as desired.

Note that by assuming that all $S_i$ are 4-sets,
we ended up with all of them sharing the two elements 2 and 3.
If we instead assumed that
$S_i$ were $n$-sets with $n \ge 5$,
the proof would remain the same
except that we would end up with
all $S_i$ sharing more than two elements.
Hence there is no loss of generality
in assuming that $S_i$ are 4-sets
(in fact, this shows that
we could even assume they are 2-sets).

To prove that $\Delta_4$ is not JIS,
we start with the same assumptions that
$S_1$, $S_2$, and $S_3$ share an \nsub,
$S_2$, $S_3$, and $S_4$ share an \nsup,
and
$S_1 = 1234$,
$S_2 = 1235$,
$S_3 = 1236$, and
$S_4 = 2356$.
Now, $S_3$, $S_4$, and $S_5$ must share an \nsub.
So $S_5$ must contain $S_3 \cap S_4 = 236$.
Since $v_5$ is adjacent to neither $v_1$ nor $v_2$,
$S_5$ can contain neither 1 nor 4 nor 5.
Hence, without loss of generality,
$S_5 = 2367$.
Continuing,
$S_4$, $S_5$, and $S_6$ must share an \nsup.
So $S_6$ must be a 4-subset of
$S_4 \cup S_5 = 23567$;
i.e., we must drop one element from 23567 to get $S_6$.
Dropping 5 or 7 gets us back to $S_4$ and $S_5$;
hence we must drop 2, 3, or 6.
The roles of 2 and 3 have been identical so far;
so, without loss of generality,
we must drop 2 or 6;
so $S_6 = 2357$ or $3567$.
The former is not possible since
$v_6$ and $v_2$ are not adjacent.
And the latter is ruled out by noticing
that 3567 differs from $S_1 =1234$ by three elements,
which contradicts the fact that
$v_6$ and $v_1$ have distance
two\footnote{Note that
$\Delta_4 - v_7$ \emph{is} JIS,
with $S_1$ and $S_6$ differing in three elements.
We will refer back to this point
at the very end of this section.}.
Thus we have reached a contradiction, as desired.

\begin{figure}[ht]

 \centering
 \includegraphics[width=60mm]{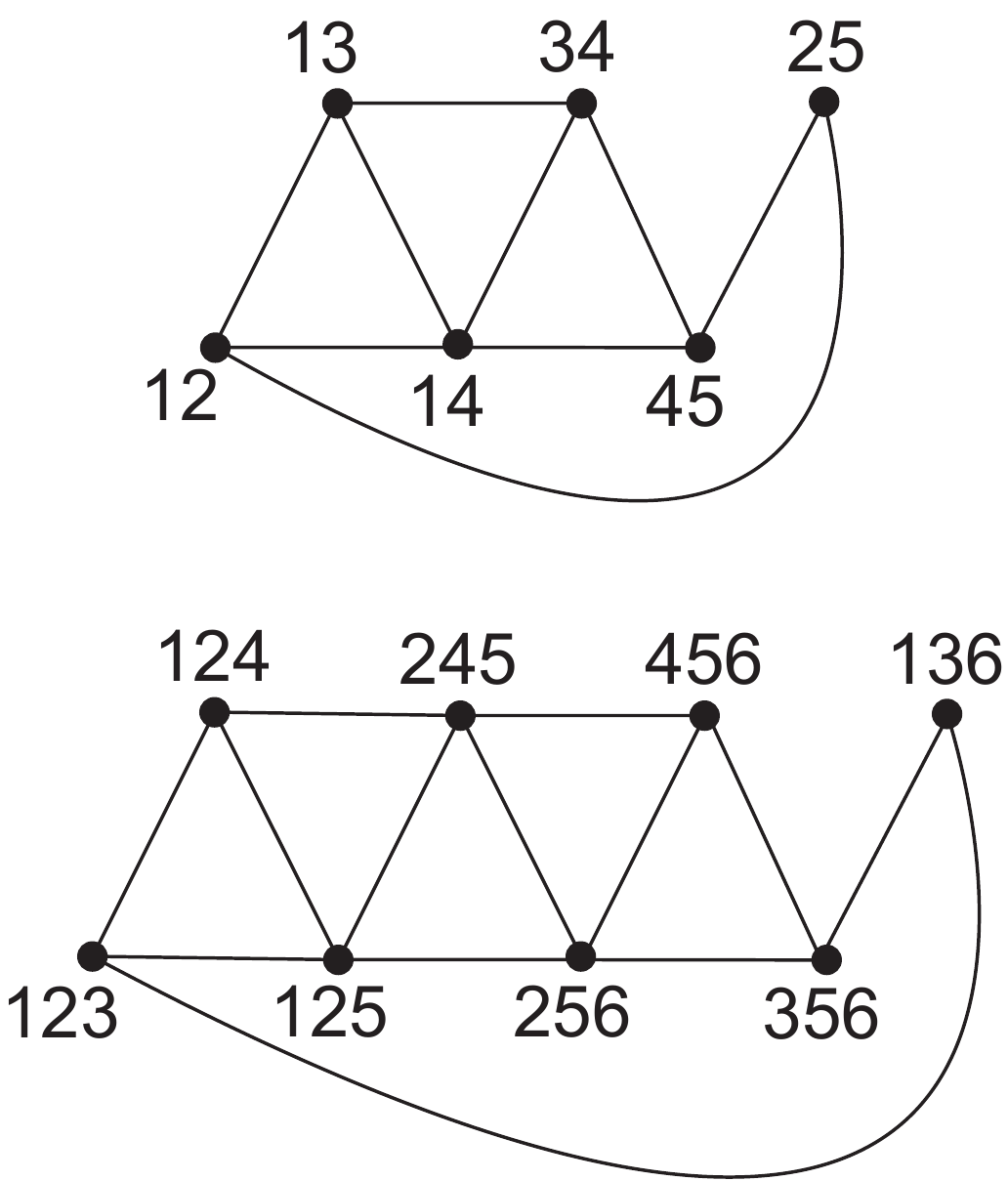}
 \caption{$\Delta_3$ and $\Delta_5$ realized as JIS graphs.
 \label{Delta35AreJISFigure}
 }
\end{figure}

Part (ii) of the proposition is proved in
\figref{Delta35AreJISFigure},
which
shows sets that realize
$\Delta_3$ and $\Delta_5$
as JIS graphs.
For the sake of compactness,
braces and commas are omitted from the sets.
}

We end this section with the following definition and question.
Let $G$ be a JIS graph,
and suppose $F = \{S_u \st u \in V(G)\}$
realizes $G$ as a JIS graph.
We define the \dfn{$F$-distance} between
two vertices $v$ and $w$ of $G$
to be $d_F(v,w)=|S_v \setsubtr S_w|$.
It is easy to show this distance function
is indeed a metric.
The \dfn{JIS-diameter} of $G$
is defined as
$$\max_{v,w \in V(G)} \min_F \{d_F(v,w)\}$$
where the minimum is taken over all families $F$ that realize
$G$ as a JIS graph.

\xquote{
\emph{Question:}
Do there exist JIS graphs with arbitrarily large JIS-diameter?
}

From the proof of Proposition~\ref{Delta2and4NotJISTheorem}
and the footnote in it,
it follows that
$\Delta_4$ minus the degree-2 vertex $v_7$
has JIS-diameter 3:
$S_1 = 1234$,
$S_2 = 1235$,
$S_3 = 1236$,
$S_4 = 2356$,
$S_5 = 2367$, and
$S_6 = 3567$,
i.e., $v_1$ and $v_6$
have $F$-distance 3.

\section{An algorithm for recognizing JIS graphs}

As mentioned in the introduction,
the following theorem provides for
an algorithm for deciding if a graph is JIS
by doing a bounded exhaustive search.

\thm{ \label{UpperboundTheorem}
\UpperboundTheoremStatement
}

\begin{proof}{
Let $G$ be a JIS graph of order $n$
with $c$ connected components.

\medskip
\emph{Case 1.}
Assume $c = 1$,
i.e., $G$ is connected.
In this case we will prove a slightly stronger result,
which we will use in the proof of Case~2:
\begin{quote}
Claim: $G$ is isomorphic, for some $m \le n$, to an induced subgraph of $J(m,2n-1)$.
\end{quote}
The case $n=1$ is trivial;
so we assume $n \ge 2$.
Since $G$ is connected,
there exists an ordering
$v_1, v_2, \cdots, v_n$ of the vertices of $G$
such that for each $i \ge 2$, $v_i$ is adjacent
to at least one of $v_1, \cdots, v_{i-1}$.
Since $G$ is JIS,
for some $k \geq 1$
there exist $k$-sets $\{S_{1}, \cdots, S_{n}\}$
that realize $G$ as a JIS graph,
where $S_i$ corresponds to the vertex $v_{i}$.
Since $v_1$ and $v_2$ are adjacent,
$| S_{1} \cap S_{2} | = k-1$.
Since $v_3$ is adjacent to at least one of $v_1$ and $v_2$,
$| S_{1} \cap S_{2} \cap S_{3} | \ge k-2$.
Continuing this way,
we see that
$| S_{1} \cap \cdots \cap S_{n} | \ge k-(n-1)$.
Let $S'_{i} = S_{i} \setsubtr (S_{1} \cap \cdots \cap S_{n})$ for $1 \le i \le n$.
Then for all $i$, $|S'_{i}| =m$ where $m \le k -(k-(n-1)) = n-1$,
and it is easily verified that
the family of sets $\{S'_{1}, \cdots, S'_{n}\}$
realizes $G$ as a JIS graph.

Now, since $v_1$ and $v_2$ are adjacent,
$| S'_{1} \cup S'_{2} | = m+1$.
Since $v_3$ is adjacent to at least one of $v_1$ and $v_2$,
$| S'_{1} \cup S'_{2} \cup S'_{3} | \le m+2$.
Continuing this way,
we see that
$ | S'_{1} \cup \cdots \cup S'_{n} | \le  m+n-1 \le 2n-2$,
which implies $G$ is an induced subgraph of $J(m,2n-1)$, $m \le n-1$.
(Note: we proved the inequalities
$ | S'_{1} \cup \cdots \cup S'_{n} | \le 2n-2$ and $m \le n-1$
only for $n \geq 2$, not for $n=1$.)

\medskip
\emph{Case 2.}
Assume $c \ge 2$.
Let $n_i$ be the order of the $i$th component of $G$.
Then, by Case~1 above,
for each $i$ there is a family $F_i$ of $m_{i}$-sets, $m_i \le n_i$,
that realizes the $i$th component of $G$ as a JIS graph,
such that
the union of the sets in $F_i$
has at most $2n_i -1 $ elements.
Thus $\bigcup F_i$ has at most
$2n-c$ elements.

We can assume $m_1 \ge m_i$ for all $i$.
We can also assume that
for all $i \ne j$,
every set in the family $F_{i}$ is disjoint from every set in $F_{j}$.
To make all sets in all the families have the same size,
for each $i$ such that $m_1 > m_i$
we add the first $m_1 - m_i$ elements
of the first set in $F_1$
to every set in $F_i$.
After adding these extra elements,
we must make sure that
sets corresponding to vertices in different components of $G$
do not share immediate subsets.
This will automatically be true
for sets that had two or more elements
before the extra elements were added,
but not for singletons.
We remedy this
by adding, for each $i$,
an element $e_i$
to every set in $F_i$,
where $e_1, \cdots, e_c$
are distinct elements not already in any set in any $F_i$.
Let $F = \bigcup F_i$.
Then $G$ is realized as a JIS graph by $F$,
which is a family of
$(m_1+1)$-sets whose union has at most $2n-c + c = 2n$ elements,
where  $m_1 + 1 \le n_1 +1 \le n$.
Thus $G$ is an induced subgraph of $J(m, 2n)$
where $m = m_1+1 \le n$.

}\end{proof}

Remark.
It is not difficult to modify the above proof in Case~1 to show that
if $G$ is connected,
then it is an induced subgraph of $J(n,2n)$.
It would be interesting to see for which graphs the bounds $n$ and $2n$ can be lowered.
Note that if $G$ consists of exactly $n \geq 2$ vertices of degree zero,
then the bound $2n$ is optimal.

\section{Edge move distance graphs and JIS graphs}

Since the 1970's
many authors have written on various metrics
defined on sets of graphs
(e.g., see \cite{BGMW}, \cite{ChGHJ}, \cite{ChGHLSW}, \cite{DD}, \cite{JohnsonM}, \cite{Kaden}, \cite{Zelinka2}).
Among them are
edge move, edge rotation, edge jump, and edge slide distances,
to which we add a new one,
edge skip distance, which we'll define later in this section.
In general,
given a metric $d$ on a set of graphs
$S= \{G_1, \cdots, G_k \}$,
the \dfn{distance graph} of $S$, denoted $D_d(S)$,
has $S$ as its vertex set,
where two vertices $G_i$ and $G_j$
are adjacent if $d(G_i, G_j) = 1$.
We will see shortly that
distance graphs associated with the
edge move metric are closely related to JIS graphs.

An \dfn{edge move} on a graph $G$
consists of removing one edge from and adding a new edge to $G$,
without changing its vertex set $V(G)$;
i.e., one edge is ``moved to a new position.''
The \dfn{edge move distance} $d_m(G,H)$
between two graphs $G$ and $H$
is defined as the fewest number of edge moves
necessary to transform $G$ into $H$,
up to isomorphism.
Note that for $d_m(G,H)$ to be defined,
$G$ and $H$ must have
the same order and the same size.
It is easy to verify that $d_m$ is a metric
on any set of graphs of given order and size.
Given a set $S$ of graphs of the same order and size,
the \dfn{edge move distance graph} of $S$, $D_m(S)$,
is the graph whose vertices are the elements of $S$,
where two vertices are adjacent
if their edge move distance is one.
When we say
a graph is an edge move distance graph
we mean it is isomorphic to one.

The connection between
JIS graphs and edge move distance graphs
can be seen by focusing on edge sets.
Let $G$ and $H$ be graphs of the same order and size,
with $n$ edges each.
If the edge sets $E(G)$ and $E(H)$
share exactly $n-1$ elements,
then $G$ and $H$ have edge move distance one.
Conversely, if $G$ and $H$ have edge move distance one,
then their vertices can be labeled such that
$E(G)$ and $E(H)$ share exactly $n-1$ elements.
At first glance,
this might seem to suggest that
a graph is JIS iff
it is isomorphic to an
edge move distance graphs.
We will show, however,
that only half (one direction) of this statement is true.

\prop{ \label{JISimpliesEMDGTheorem}
\JISimpliesEMDGTheoremStatement
}

\pf{
Let $G$ be realized as a JIS graph by
a family of $n$-sets $\{S_v \st v \in V(G) \}$.
We will construct a graph $G_v$ for each $v \in V(G)$
such that $d_m(G_v, G_w) = 1$ iff
$S_v$ and $S_w$ share an \nsub.

We can assume that
each $S_v$ consists of positive integers.
Let $k = 1 + \max\{i\in S_v \st v \in V(G) \}$,
and let $P$ be a path of length $2k$.
Denote the vertices of $P$ by
$p_0$, $p_1$, $\cdots$, $p_{2k}$.
For each $v \in V(G)$, we let $G_v$ be the graph
consisting of $P$
plus the edges $p_i p_{2k-i}$ for all $i \in S_v$.
Then it is easily verified that
for $v \ne w$, $G_v$ is not isomorphic to $G_w$,
and
$d_m(G_v, G_w) = 1$ iff
$S_v$ and $S_w$ share an \nsub.
Therefore $G$ is isomorphic to
the edge move distance graph
$D_m( \{ G_v \st v \in V(G) \} )$.

}

The converse of the above result is not true.
The reason is that
the number of edges shared by the edge sets of two graphs
depends on how their vertices are labeled,
whereas edge move distance is measured up to graph isomorphism.

\prop{
\KnMinusOneEdgeIsNotJISTheoremStatement
}

\begin{figure}[ht]

 \centering
 \includegraphics[width=80mm]{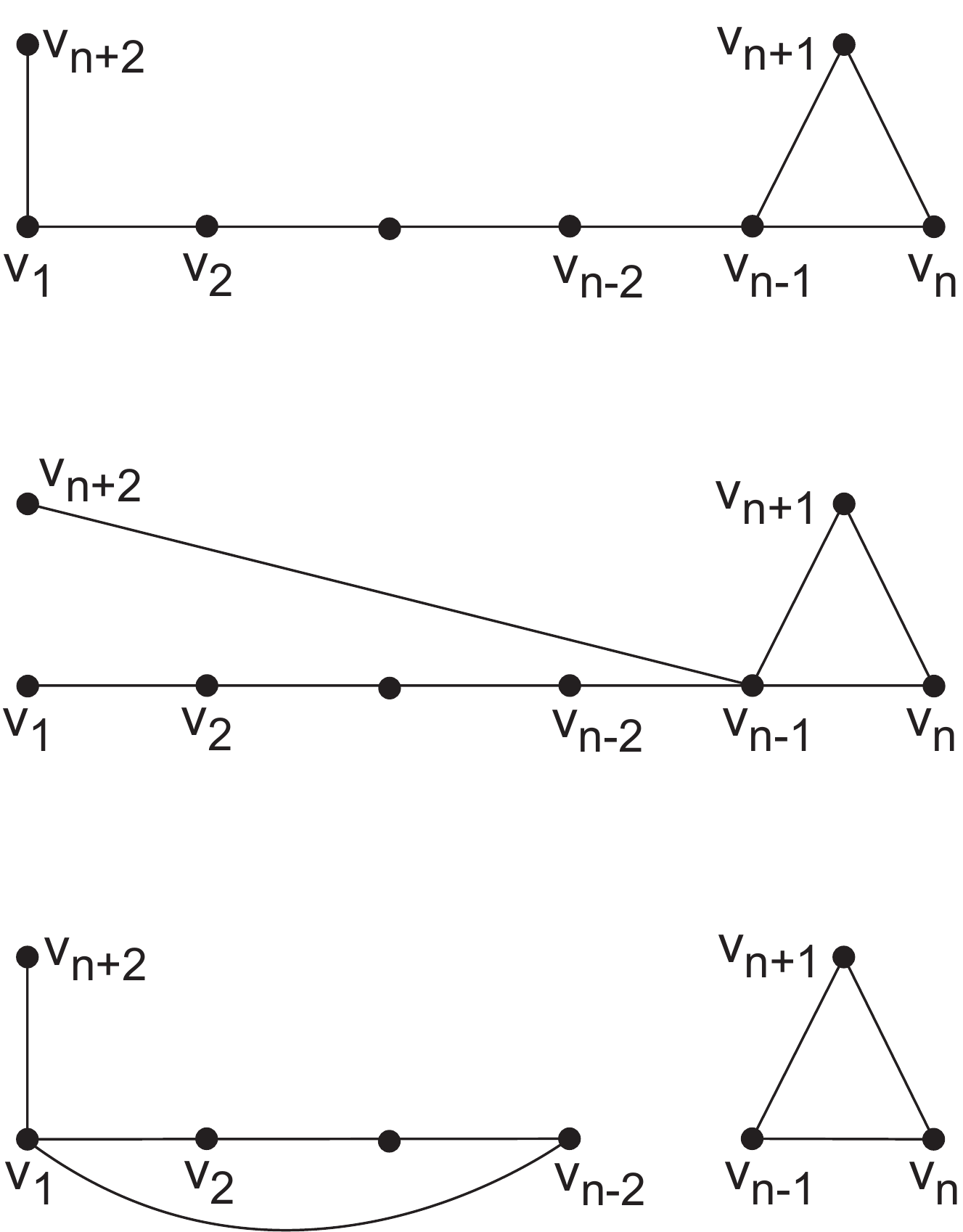}
 \caption{$Q_1$ (top), $Q_{n-1}$ (middle), and $Q_n$ (bottom)
  for $n=6$.
 \label{QiFigure}
 }
\end{figure}

\pf{
Fix $n \ge 5$,
and let $H$ be the graph obtained by removing one edge from $K_n$.
Then $H$ contains two maximal $(n-1)$-cliques
which share $n-2$ vertices.
Hence, by Part~1 of Proposition~\ref{CliquesInJISTheorem},
$H$ is not JIS.

To show that $H$ is an edge move distance graph, we construct a set
of graphs $S = \{ Q_1,  Q_2, \cdots, Q_n \}$ such that $H
\isomorphic D_m(S)$. For $1 \le i \le n$, $Q_i$ has $n+2$ vertices:
$V(Q_i) = \{v_1, v_2, \cdots, v_{n+2} \}$. For $1 \le i \le n-1$,
$E(Q_i) = \{v_k v_{k+1} \st 1 \le k \le n\} \cup \{v_{n-1}v_{n+1},
v_i v_{n+2} \}$; and $E(Q_{n}) = (E(Q_1) \cup  \{v_1 v_{n-2} \})
\setsubtr \{ v_{n-2} v_{n-1} \}$. \figref{QiFigure} shows $Q_1$,
$Q_{n-1}$, and $Q_n$ for $n=6$.

Then one readily verifies
for all $i \ne j$
except when $\{i,j\} = \{n-1, n \}$
that $Q_i$ and $Q_j$
have edge move distance one.
Thus $H$ is an edge move distance graph.

}

\subsection*{Acknowledgments}
We thank Terry A.\ McKee of Wright State University
for bringing to our attention that
the graphs we were studying
are related to Johnson graphs.
The first named author thanks
Caltech for its hospitality
while he did part of this work
there during his sabbatical leave.
The second named author thanks
the Undergraduate Research Center of Occidental College
for providing support to do this work as part of an
undergraduate summer research project.
We also thank the (anonymous) referee
for helpful suggestions.

{\scriptsize

}

}
\begin{thebibliography}{abc}

\bibitem{BGMW}
Gerhard Benad\'{e}, Wayne Goddard, Terry A.\ McKee, Paul A.\ Winter,
\emph{On Distances Between Isomorphism Classes of Graphs,}
Mathematica Bohemica \textbf{116} (1991), no.~2, 160-169.


\bibitem{ChGHJ}
Gary Chartrand, Heather Galvas, H\'{e}ctor Hevia, Mark A.\ Johnson,
\emph{Rotation and Jump Distances Between Graphs,}
Discussiones Mathematicae, Graph Theory \textbf{17} (1997) 285-300


\bibitem{ChGHLSW}
Gary Chartrand, Wayne Goddard, M.\ A.\ Henning, Linda Lesniak, Henda C.\ Swart, Curtiss E.\ Wall,
\emph{Which Graphs Are Distance Graphs?}
Ars Combin.\ \textbf{29A} (1990) 225-232


\bibitem{DD}
Michel Marie Deza, Elena Deza,
emph{Encyclopedia of Distances},
Springer-Verlag, Berlin Heidelberg, 2009.

%\bibitem{egp}
%Paul Erd\"{o}s, A.\ W.\ Goodman, Lajos P\'{o}sa,
%\emph{The representation of a graph by set intersection,}
%Canad.\ J.\ Math.\ \textbf{18} (1966) 106-112.

\bibitem{JohnsonM}
Mark A.\ Johnson,
\emph{An Ordering of Some Metrics Defined on the Space of Graphs,}
Czechoslovak Mathematical Journal \textbf{37} (1987) 75-85.



\bibitem{Kaden}
Frieder Kaden,
\emph{Graph metrics and distance-graphs,}
Graphs and other combinatorial topics (Prague, 1982),
145-158.


\bibitem{MKMM}
Terry A.\ McKee,
F.\ R.\ McMorris,
\emph{Topics in Intersection Graph Theory,}
SIAM Monographs on Discrete Mathematics and Applications, 1999.

\bibitem{Terwil}
Paul Terwilliger,
\emph{The Johnson graph $J(d,r)$ is unique if $(d,r)\not=(2,8)$.}
Discrete Mathematics \textbf{58} (1986) 175-189


%%%
%%%\bibitem{Zelinka}
%%%Bohdan Zelinka,
%%%\emph{On a certain distance between isomorphism classes of graphs,}
%%%\u{C}asopis.\ P\u{e}st.\ Mat.\ \textbf{100}(1975) 371-373.

\bibitem{Zelinka2}
Bohdan Zelinka,
\emph{Comparison of various distances between isomorphism classes of graphs},
\u{C}asopis.\ P\u{e}st.\ Mat.\ \textbf{110} (1985), no.~3, 289-293, 315.



\end{thebibliography}
\end{document}